\documentclass[11pt]{article}
\usepackage{makeidx}
\usepackage[centertags]{amsmath}
\usepackage{amsthm}
\usepackage{newlfont}
\usepackage{amsmath}
\usepackage{amsfonts,amsmath,latexsym}
\usepackage{mathrsfs}
\usepackage{amssymb}
\usepackage{amsbsy}
\usepackage{indentfirst}
\usepackage{mathrsfs}
\usepackage{graphicx}
\usepackage{epstopdf}
\usepackage{subfigure}
\usepackage{caption}
\usepackage{grffile}
\usepackage{float}
\usepackage{amssymb, amsmath}
\usepackage{xcolor}
\usepackage{enumitem}
\usepackage{CJK, CJKnumb, CJKulem, verbatim}
\usepackage[colorlinks, citecolor=blue, linkcolor=blue]{hyperref}
\usepackage{diagbox}
\usepackage[T1]{fontenc}

\titlepage
\newlength{\defbaselineskip}
\newcommand{\ep}{\epsilon}

\newcommand{\be}{\begin{eqnarray}}
\newcommand{\ee}{\end{eqnarray}}
\newcommand{\bestar}{\begin{eqnarray*}}
\newcommand{\eestar}{\end{eqnarray*}}

\newcommand{\ignore}[1]{}
{} \theoremstyle{plain}
\newtheorem{thm}{Theorem}[section]

\newtheorem{lemma}{Lemma}[section]

\theoremstyle{definition}

\newtheorem{rem}{Remark}[section]
\numberwithin{equation}{section}

\def\>{\geq}
\def\<{\leq}

\allowdisplaybreaks
\graphicspath{{./graphics/}}

\begin{document}

\title{The elephant random walk with  
proportionally growing recent memory}
 \author{Zhishui Hu
\thanks{E-mail addresses: huzs@ustc.edu.cn} \qquad  Shuwen Xue\\
\\
Department of Statistics
of Finance, School of Management\\
 University of Science and Technology of
China\\
 Hefei, Anhui 230026, China
}
\date{}
\maketitle

\begin{abstract} 

In this work, we study an elephant random walk with recent restricted memory, where at step $n$ the
elephant remembers only the last $m_n=\max\{\lfloor cn\rfloor,1\}$ steps
for $c\in(0,1)$.  
We show that the random walk undergoes a phase transition at the critical point $p_c=(3+\sqrt{1-c})/4$, which depends on the memory proportion $c$.
For $0<p<p_c$, $S_n/\sqrt{n}$ is asymptotically normal; for $p=p_c$, $S_n/\sqrt{n\log n}$ is asymptotically normal; and for $p_c<p<1$, $S_n/n^{\rho}$  converges almost surely and in $L^2$ to a random variable for some constant $\rho\in (1/2, 1)$. This phase transition differs from that of the classical elephant random walk and its remote-memory counterpart, both of which have the same fixed critical point $3/4$.

 \vskip 0.2cm


\end{abstract}

\section{Introduction and main result}

The  elephant random walk (ERW)  is a discrete-time random walk on $\mathbb{Z}$ with long-range memory, introduced by Sch\"{u}tz and Trimper \cite{Schutz2004}.  
Let $X_1$  be  uniformly distributed on $\{-1,1\}$, and,
for each step $n\ge 2$,  define recursively
\be 
X_n=\left\{
\begin{array}{ll}
+X_{U_n}, & \text{with probability} ~p,\\
-X_{U_n}, & \text{with probability} ~1-p,
\end{array}
\right.  \label{defrecur}
\ee 
where $\{U_n, n\ge 2\}$ are independent and $U_n$ is  uniform on  the integers $1,2,\cdots, n-1$ for each $n\ge 2$.
Then $S_n=\sum_{k=1}^n X_k$ is an ERW,  and the parameter $p$ is called the reinforcement parameter. 

The limiting behavior of the ERW has been extensively studied;
see, e.g., \cite{Baur2016, Bercu2018, Bercu2021, Coletti2017a, Coletti2017b, Kubota2019, Kursten2016, Qin2025}. 
In particular,   Bercu \cite{Bercu2018} established the following limit theorem.

\begin{thm}\label{thBercu}  \cite{Bercu2018} Let $S_n$ be an ERW with reinforcement parameter $p$.  
\begin{enumerate}
\item[$(i)$] (Diffusive regime)  If $0<p<3/4$,  then
  $$\frac{S_n}{\sqrt{n}}\xrightarrow[\quad]{d}N\Big(0, \frac{1}{3-4p}\Big);$$ 

\item[$(ii)$] (Critical regime) If $p=3/4$, then  $$\frac{S_n}{\sqrt{n\log n}}\xrightarrow[\quad]{d} N(0,1);$$  

\item[$(iii)$] (Superdiffusive regime) If $3/4<p<1$, then 
\bestar
\frac{S_n}{n^{2p-1}}\to L ~~\text{a.s. and in}~ L^2, 
\eestar
where $L$ is a non-Gaussian random variable.

\end{enumerate}
\end{thm}

In the ERW, the  elephant remembers all  of its past steps: at each step, it samples uniformly from all previous steps. Theorem \ref{thBercu} shows that long-range memory affects the dynamics of the random walk, leading to a crossover from diffusive to superdiffusive behavior as $p$ increases. In this work, we focus on the case in which the elephant's memory is restricted to a fraction of the past, and investigate how this restriction modifies the scaling properties of the random walk.

  In the  ERW with restricted memory, $U_n$ is uniform on $\mathfrak{M}_n$ in \eqref{defrecur}, where $\mathfrak{M}_n\subset\{1,2,\cdots,n-1\}$ is the elephant's memory over the first $n-1$ steps. This model has been explored by simulation in \cite{ Cressoni2007, DaSilva2006}, where a given fraction of the remote/recent past is remembered.
The theoretical properties  have also attracted considerable attention recently, including  \cite{Aguech2024, BenAri2021,   Gut2021, Gut2022, Gut2023, Gut2026, Nakano2025, Roy2026, Vazquez2025},  among others, covering the cases in which the elephant remembering only the remote past, only  the recent past, or a mixture of both.

For the ERW remembering the remote past, i.e., $\mathfrak{M}_n=\{1,2,\cdots,m_n\}$ with $1\le m_n\le n-1$, there is no phase transition when $m_n=m~(n\ge m+1)$ is fixed (see \cite{Gut2021}), whereas when $m_n$ is increasing and $m_n/n\to c\in[0,1]$, a phase transition occurs at $p=3/4$, analogous to Theorem \ref{thBercu} (see \cite{Aguech2024, Gut2022}).

For the ERW remembering the recent past, i.e., $\mathfrak{M}_n=\{n-m_n, n-m_n+1, \cdots, n-1\}$ with $1\le m_n\le n-1$,  
 $S_n/\sqrt{n}$ is asymptotically normal when $m_n=m~(n\ge m+1)$ is fixed (see  \cite{BenAri2021, Gut2021}).
Moreover, Gut and  Stadtm\"uller \cite{Gut2022} conjectured that asymptotic normality still holds for  increasing $m_n$. 

The conjecture in \cite{Gut2022} suggests that the two scenarios of the ERW remembering the remote/recent past differ significantly, especially when $m_n \to \infty$. Indeed, when $m_n$ is proportional to $n$,  the simulation results in \cite{Cressoni2007} and \cite{DaSilva2006}  also indicate a  difference between the two cases. However,  no theoretical result is available for $m_n \rightarrow \infty$ in the recent-memory case.  The aim of this paper is to address this problem. We assume that $
m_n=\max\{\lfloor cn\rfloor, 1\}
$ for  $c\in [0,1]$. Since the cases $c=0$ and $c=1$ have been studied in \cite{Gut2021}, we restrict our attention to $c\in(0,1)$.

To state our main result, let
 $w(t)$ be the fundamental solution of  the following linear retarded differential difference equation
\be
w'(t)=\frac{2p-1}{c}(w(t)-w(t+\log (1-c))),~~t\ge 0,  \label{RDD0}
\ee
that is,  $w(t)$ solves \eqref{RDD0} and satisfies the initial condition
\be 
w(t)=\left\{
\begin{array}{ll}
0, & t<0,\\
1, & t=0.
\end{array}   
\right.   \label{initialcond}
\ee 
For further properties of $w(t)$, 
we refer to \cite{Hale1993}  and to Section \ref{sect3} below.

Our main result is as follows.

\begin{thm} \label{th1}
Let $S_n$ be an ERW with restricted memory $\mathfrak{M}_n=\{n-m_n, n-m_n+1, \cdots, n-1\}$,  where $m_n=\max\{\lfloor cn\rfloor, 1\}$ for some $c\in (0,1)$.
Define 
\be
p_c=\frac{3+\sqrt{1-c}}{4}.  \label{pc}
\ee
\begin{enumerate}
\item[$(1)$] (Diffusive regime) If $0<p<p_c$, then
\be 
\frac{S_n}{\sqrt{n}}\xrightarrow[\quad]{d} N(0, \sigma^2(c,p)), 
\label{CLTsubcritical}
\ee 
where \bestar 
\sigma^2(c, p) = \int_{0}^{\infty} \frac{w^2(u)}{e^u} du\in (0,\infty). 
\eestar

\item[$(2)$] (Critical regime) If $p=p_c$, then
\be 
\frac{S_n}{\sqrt{n\log n}}\xrightarrow[\quad]{d} N(0, \sigma^2_c), 
\label{CLTcritical}
\ee 
where 
\bestar 
\sigma^2_c=\Big(1+\frac{(1-c+\sqrt{1-c})\log (1-c)}{2c}\Big)^{-2}.
\eestar 

\item[$(3)$] (Superdiffusive regime) If $p_c<p<1$, then
there exists a random variable $W$ such that
\be 
\frac{S_n}{n^{\rho}} \longrightarrow W\qquad \mbox{a.s. and in}~L^2. \label{superconv}
\ee 
where $\rho \in (1/2, 1)$ is the unique positive root of the equation
\bestar
z=\frac{2p-1}{c} (1-(1-c)^z).
\eestar
\end{enumerate}
\end{thm}

\begin{rem}
Theorem \ref{th1} shows that $S_n$ exhibits a phase transition at $p_c$.
In contrast to the classical ERW and the ERW with the remote past, whose common critical point is fixed at  $3/4$, 
 the critical point  $p_c$ here is decreasing in $c$ and tends to $3/4$ as $c\to 1$. On the other hand,   $p_c\to 1$ as $c\to 0$. This coincides with  the result established in \cite{Gut2021} for $c=0$ ($m_n=n-1$), which shows that $S_n/\sqrt{n}$ is asymptotically normal for
all $p\in (0,1)$, and hence  no phase transition occurs in this case. In addition, Theorem \ref{th1} also shows that the conjecture in \cite{Gut2022} is false.
\end{rem}

The main difficulty in proving our main theorem is that  we cannot construct a martingale directly as in \cite{Bercu2018}. Instead, we adopt a different martingale  approach and express  $S_n$ as $\sum_{k=1}^n w_{n,k}\ep_k$,  where $\{w_{n,k}\}$ is   a triangular array of weights and $\{\ep_k\}$ is a martingale difference sequence.  To investigate the asymptotic behavior of these weights,  we introduce suitable interpolations of  $\{w_{n,k}\}$  and establish their uniform convergence on any compact interval. The limiting function  $w(t)$ is the fundamental solution of the retarded differential difference equation \eqref{RDD0}.  The details of this  approach and the proof  strategy are  described in Section \ref{sect2}. In Section \ref{sect3}, we establish some basic properties of  $\{w_{n,k}\}$ and  $w(t)$. These properties are then used to prove our main theorem in Section  \ref{sect4}.

Throughout this paper, we adopt the following notation. 
For an event $E$, let $I(E)$ denote its indicator function. 
For a   real sequence $a_n$ and a  positive  sequence  $b_n$, we write
$a_n \sim b_n$ if $a_n/b_n \to 1$,
$a_n = o(b_n)$ if $a_n/b_n \to 0$,
$a_n = O(b_n)$ if $a_n/b_n$ is bounded,
as $n \to \infty$. The same notation applies to quantities depending on multiple indices or a real variable, with limits taken accordingly.
We use  $\xrightarrow{d}$ and $\xrightarrow{L^2}$  to denote convergence in distribution and in mean square, respectively.
Finally,  $C$ denotes a generic positive constant independent of 
$n$ and $k$, which may take different values at different occurrences.

\section{Martingale  approach and proof strategy} \label{sect2}

For notational simplicity, in what follows we write
\bestar 
  a:= 2p-1, \qquad L:= -\log(1-c), \qquad \lambda:= aL/c.
\eestar

The martingale approach has been widely employed in the study of the ERW.  Let $\mathscr{F}_k=\sigma(X_1, \cdots, X_k)$ for $k \ge 1$.   Then,  for the classical ERW,  
\bestar
\mathbb{E}(X_{n}|\mathscr{F}_{n-1})=  \frac{aS_{n-1}}{n-1},\qquad n\ge 2,
\eestar
and consequently,   $\{M_n:=a_n S_n, n\ge 1\}$ is a martingale (see \cite{Bercu2018}
for more details), where
\bestar
a_n=\frac{\Gamma(n)\Gamma(2p)}{\Gamma(n+2p-1)}.
\eestar
 However,  for the ERW with restricted memory $\mathfrak{M}_n$,
we have
\bestar
\mathbb{E}(X_{n}|\mathscr{F}_{n-1})= \frac{a\sum_{i\in \mathfrak{M}_n}X_i}{|\mathfrak{M}_n|}, \qquad n\ge 2.
\eestar
In this case, unlike the classical ERW, a martingale cannot be constructed directly through a deterministic normalization of $S_n$.

In this paper, we adopt a different martingale  technique. 
Our approach relies on the fact that $S_n$ can be written as a sum of martingale differences with respect to the natural filtration generated by $\{X_k\}$.
More precisely,   for any $n\ge 1$,
\be
S_n=\sum_{k=1}^n (\mathbb{E}(S_n|\mathscr{F}_{k})- \mathbb{E}(S_n|\mathscr{F}_{k-1})),  \label{MD}
\ee
where   $\mathscr{F}_0=\{\emptyset, \Omega\}$, and
$\{\mathbb{E}(S_n|\mathscr{F}_{k})- \mathbb{E}(S_n|\mathscr{F}_{k-1}), 1\le k\le n\}$ is a martingale difference sequence with respect to $\{\mathscr F_k, 1\le k\le n\}$.
For an introduction to this method and its  applications to concentration inequalities, we refer to Section 2.2 in \cite{Wainwright2019}.

To facilitate the use of martingale theory,    Lemma  \ref{lemma1} below provides a further representation
of $S_n$ in terms of the innovations 
\be
\ep_k:=X_k-\mathbb{E}(X_k|\mathscr{F}_{k-1}),\qquad k\ge 1.  \label{innovation}
\ee
The desired representation can be established by  using \eqref{MD}
and verifying that 
$\mathbb{E}(S_n|\mathscr{F}_{k})- \mathbb{E}(S_n|\mathscr{F}_{k-1})
=w_{n,k}\ep_k$. Here we present a simpler proof based on a recursive argument.

\begin{lemma} \label{lemma1}
Assume that the weights $\{w_{n,k}, k\ge 1, n\ge 1\}$  satisfy that $w_{n,k}=0$ for $k>n$,
$w_{n,n}=1$ and for $k<n$,
\be  
w_{n,k}=w_{n-1,k}+\frac{a}{m_n}(w_{n-1, k}-w_{n-m_n-1, k}).  \label{windu}
\ee 
Then, for any $n\ge 1$, 
\be 
S_n=\sum_{k=1}^n w_{n,k} \ep_k,   \label{Mar0}
\ee 
where  $\varepsilon_k$ is defined in \eqref{innovation}.
\end{lemma}

\begin{proof} 
Recalling the definition of $X_n$ and using \eqref{innovation}, we obtain that for $n\ge 2$,
\be 
\ep_n=X_n-\mathbb{E}(X_{n}|\mathscr{F}_{n-1})=X_n-\frac{a}{m_n}\sum_{k=n-m_n}^{n-1} X_k=X_n-\frac{a}{m_n}(S_{n-1}-S_{n-m_n-1}).  \label{Mar1}
\ee

We now prove \eqref{Mar0} by induction on $n$. The result holds when $n=1$.
Assume that \eqref{Mar0} holds for $n\le \ell-1$, where $\ell\ge 2$. Then, by
\eqref{windu} and \eqref{Mar1},  we have
\bestar
\sum_{k=1}^{\ell} w_{\ell,k} \ep_k&=&\ep_{\ell}+\sum_{k=1}^{{\ell}-1} w_{{\ell}-1,k} \ep_k+\frac{a}{m_{\ell}}\Big(\sum_{k=1}^{{\ell}-1} w_{{\ell}-1, k}\ep_k-\sum_{k=1}^{{\ell}-1} w_{{\ell}-m_{\ell}-1, k}\ep_k\Big)\\
&=& \ep_{\ell}+S_{{\ell}-1}+ \frac{a}{m_{\ell}}(S_{{\ell}-1}-S_{{\ell}-m_{\ell}-1})=S_{\ell}.
\eestar
Thus \eqref{Mar0} holds for $n=\ell$, completing the induction.
\end{proof}

In what follows, we briefly outline the strategy for proving the main result of this paper.

In the proof, we apply martingale limit theorems to the array $\{w_{n,k} \ep_k,  n\ge 1, 1\le k\le n\}$. This requires us to investigate the asymptotic behavior of $w_{n,k}$  as  $n\to\infty$. Although $w_{n,k}$ is defined explicitly through \eqref{windu}, not all the required asymptotic properties can be obtained directly from this recursion. Instead, we consider the sequence $\{w_{n,k}, n\ge k\}$ as a whole and construct a continuous interpolation on $[0,\infty)$.

Let $\tau_{n,k}=\log (n/k)$ for $n\ge k$.
we define $w_k\in C[0,\infty)$ such that $w_k(\tau_{n,k})=w_{n,k}$ and $w_k$ is linear on each interval $[\tau_{n-1,k}, \tau_{n,k}]$. Specifically,
for $t\ge 0$,
\be 
w_k(t)=w_{\lfloor ke^t\rfloor,k}+\frac{t-\log (\lfloor ke^t\rfloor/k)}{\log (1+(\lfloor ke^t\rfloor)^{-1})}(w_{\lfloor ke^t\rfloor+1,k}-w_{\lfloor ke^t\rfloor,k}).  \label{wkt}
\ee
For convenience,  set $w_k(t)=0$ for $t<0$.

Using some  simple  properties of $w_{n,k}$, 
 Lemma \ref{lemma3.2} below shows that
for any fixed $K>0$, the family
$\{w_k(t): 0\le t\le K\}$ is relatively compact  in $(C[0,K], \rho_K)$,  where $C[0,K]$ is the space of continuous functions on $[0,K]$ 
and $\rho_K$ denotes the uniform metric on $C[0,K]$. 
We then investigate the possible subsequential limits of  $\{w_k\}$ as $k\rightarrow \infty$.

A further analysis (see Lemma \ref{Lemmawk}) shows that $w_k$ is absolutely continuous and satisfies 
\bestar
w_k'(t)=\frac{a}{c}\big(w_k(t)-w_k(t-L)\big)+R_k(t),  
\eestar
for $t\notin\{\tau_{n,k}:\, n\ge k\}$,
where $R_k(t)$ is the remainder term. This observation suggests that the limiting function should satisfy the retarded differential difference equation 
\bestar
w'(t)=\frac{a}{c}\big(w(t)-w(t-L)\big),  
\eestar
together with  the initial condition \eqref{initialcond}. By analyzing the solution of this equation and using the relationship between $w_k$ and $w$ (see \eqref{wktwt}):
\bestar
w_k(t)=w(t)+\int_0^t w(t-s)R_k(s)\,ds,  
\eestar
we establish in Lemma  \ref{lemma3.6} that $w_k\to w$ in $(C[0,K],\rho_K)$,   and derive precise estimates for $w_{n,k}$ in Lemma \ref{lemma3.7}. Finally,  these results, together with  limit theorems for martingale difference arrays,  enable us to complete the proof of the main result.  The detailed arguments are presented in Sections \ref{sect3} and \ref{sect4}.

\section{Preliminaries}  \label{sect3}

\begin{lemma} \label{lemma3.1}
Let $w_{n,k}$ be defined as in Lemma \ref{lemma1}.
\begin{enumerate}
\item[$(1)$] If $p\le 1/2$, then 
$0< w_{n,k}\le 1$ for any  $n\ge k$. 

\item[$(2)$]
If  $p>1/2$,  then $w_{n,k}\ge 1$ for any $n\ge k$, and $w_{n,k}$
is non-decreasing in $n$.  Moreover,  we have $
 w_{n, k}\le  w_{n, j}w_{j,k}
$ for any $k \le j \le n$.

\end{enumerate}
\end{lemma}

\begin{proof} 
(1)  If $p\le 1/2$, then $-a/m_n\in [0,1)$, and by  \eqref{windu},
\bestar
w_{n,k}=\Big(1+\frac{a}{m_n}\Big)w_{n-1, k}+ \frac{-a}{m_n}w_{n-m_n-1, k}.
\eestar
This means that $w_{n,k}$ 
is a convex combination of $w_{n-1, k}$ and $w_{n-m_n-1, k}$.
Since $w_{k,k}=1$ and $w_{n,k}=0$ for $n<k$, induction on $n$
gives $w_{n,k}\in (0,1]$ for all $n\ge k$.

(2)  Assume $p>1/2$. 
Then $a>0$, and by \eqref{windu},
\bestar
w_{n,k}-w_{n-1,k}=\frac{a}{m_n}(w_{n-1, k}-w_{n-m_n-1, k}).
\eestar
Since $w_{k,k}=1$ and $w_{n,k}=0$ for $n<k$,
we obtain by induction that  $w_{n,k}\ge w_{n-1,k}$.
Hence $w_{n,k}\ge w_{k,k}\ge 1$  for $n\ge k$.

 Fix $k\le j$ and define $q_{n}:=w_{n,j}/w_{n,k}$.
For $n>j$, by \eqref{windu},
\bestar
q_{n}-q_{n-1}&=&\frac{w_{n, j}w_{n-1,k}-w_{n,k}w_{n-1,j}}{w_{n,k}w_{n-1,k}}\\
&=& \frac{(w_{n, j}-w_{n-1,j})w_{n-1,k}-(w_{n,k}-w_{n-1,k})w_{n-1,j}}{w_{n,k}w_{n-1,k}}\\
&=& \frac{a}{m_n} \frac{(w_{n-1, j}-w_{n-m_n-1, j})w_{n-1,k}
- (w_{n-1, k}-w_{n-m_n-1, k})w_{n-1,j}}{w_{n,k}w_{n-1,k}}\\
&=&\frac{aw_{n-m_n-1, k}}{m_nw_{n,k}} (q_{n-1}-q_{n-m_n-1}).
\eestar
For $p>1/2$,  the factor $aw_{n-m_n-1,k}/(m_nw_{n,k})$ is nonnegative.
Also $q_{j}=1/w_{j,k}>0$, and $q_{n}=0$ for $n<j$.
By induction,   $q_n$ is non-decreasing  for all $n\ge j$.
Therefore, for $n\ge j$,
\bestar
\frac{w_{n,j}}{w_{n,k}}=q_n\ge q_j=\frac{w_{j,j}}{w_{j,k}}=\frac{1}{w_{j,k}}.
\eestar
Multiplying by $w_{n,k}w_{j,k}>0$ yields $w_{n,k}\le w_{n,j}w_{j,k}$.
\end{proof}

\begin{lemma} \label{lemma3.2}
Let $w_k$ be defined as in Section \ref{sect2}. Then for any fixed $K>0$,
$\{w_k(t): 0\le t\le K\}$ is relatively compact  on the space $(C[0,K], \rho_K)$.
\end{lemma}

\begin{proof}  We first prove the uniform  boundedness of $\{w_k\}$.

By Lemma  \ref{lemma3.1}, if $p\le 1/2$, then $0<w_{n,k}\le 1$ for
$n\ge k$. Hence, $|w_k(t)|\le 1$ for all $t\in[0,K]$.
Now assume $p>1/2$. Then $a>0$, and 
by \eqref{windu},
\bestar
w_{n,k}=\Big(1+\frac{a}{m_n}\Big)w_{n-1, k}+ \frac{-a}{m_n}w_{n-m_n-1, k}
\le \Big(1+\frac{a}{m_n}\Big)w_{n-1, k}.
\eestar
Since $m_j>0$ and $m_j\sim cj$ as $j\to \infty$,
  there exists $c_0>0$ such that
$m_j\ge c_0j$ for all $j\ge 1$.
Therefore, for $n\ge k$,
\be 
w_{n,k}&\le& \prod_{j=k+1}^{n} \Big(1+\frac{a}{m_{j}}\Big)
\le  \prod_{j=k+1}^{n} \Big(1+\frac{a}{c_0j}\Big)
\le \exp\Big\{\sum_{j=k+1}^n \frac{a}{c_0j}\Big\}\nonumber\\
&\le&    \exp\Big\{\int_{k}^n \frac{a}{c_0x}dx\Big\}
=\Big(\frac{n}{k}\Big)^{a/c_0}.  \label{wnkbound0}
\ee 
Consequently,  for  $0\le t\le K$,
\bestar 
|w_k(t)|
\le \max_{k\le n\le 2e^K k} |w_{n,k}|\le  (2e^K)^{a/c_0}.
\eestar 
Define 
\bestar
C_K=\left\{
\begin{array}{ll}
1, & p\le 1/2,\\
(2e^K)^{a/c_0}, & p>1/2.
\end{array}
\right.
\eestar
Then $\{w_k(t): 0\le t\le K\}$ is uniformly bounded by $C_K$.

We next prove the uniform equicontinuity of $\{w_k\}$.
For any  $k\le n_1<n_2\le 2e^K k$,
using\eqref{windu} and the bound
$|w_{n,k}|\le C_K$ for $k\le n\le 2e^K k$, we have  
\bestar
|w_{n_2,k}-w_{n_1, k}|&\le& \sum_{n=n_1+1}^{n_2} |w_{n,k}-w_{n-1, k}| 
 \le 2|a|C_K\sum_{n=n_1+1}^{n_2} \frac{1}{m_n}\\
&\le& \frac{2|a|C_K}{c_0}\sum_{n=n_1+1}^{n_2} \frac{1}{n}\le \frac{4|a|C_K}{c_0}\sum_{n=n_1+1}^{n_2} \int_{n}^{n+1}\frac{dx}{x}
= \frac{4|a|C_K}{c_0}\log \Big(\frac{n_2+1}{n_1+1}\Big).
\eestar
By the definition of  $w_k(t)$ as a linear interpolation of $w_{n,k}$, we have
$$|w_k(t)-w_k(s)|\le \max\{|w_{\lfloor ke^t\rfloor+i,k}-w_{\lfloor ke^s\rfloor+j,k}|: ~i,j\in \{0,1\}
\}.$$ 
Hence, for $0\le s<t\le K$ with $t-s<\delta$,
\bestar
|w_k(t)-w_k(s)|
&\le&   \frac{4|a|C_K}{c_0}\log \Big(\frac{\lfloor ke^t\rfloor+2}{\lfloor ke^s\rfloor+1}\Big)\\
&\le&  \frac{4|a|C_K}{c_0}\log \Big(\frac{ke^t+2}{ke^s}\Big)\\
&\le& \frac{4|a|C_K}{c_0}\log (e^{\delta}+2/k).
\eestar
Then, for any $\varepsilon>0$, there exist $\delta_0>0$ and $k_0\in \mathbb{N}$ such that  for $\delta<\delta_0$ and $k>k_0$,
$$\sup_{s,t\in [0,K], |t-s|<\delta}|w_k(t)-w_k(s)|<\varepsilon.$$
For the finitely many $k\le k_0$,  each $w_k(t)$ is continuous on $[0,K]$,
hence uniformly continuous. Therefore,  there exists $\delta_1>0$
such that  for $\delta<\delta_1$ and all $k\le k_0$, $$\sup_{s,t\in [0,K], |t-s|<\delta}|w_k(t)-w_k(s)|<\varepsilon.$$  Take $\delta<\min\{\delta_0,\delta_1\}$. Then
$$\sup_{k\in \mathbb{Z}}\sup_{s,t\in [0,K], |t-s|<\delta}|w_k(t)-w_k(s)|<\varepsilon.$$  This proves uniform equicontinuity.
By the Arzel\`{a}-Ascoli theorem (see  \cite[Theorem 7.2]{Billingsley1999}),
$\{w_k(t): 0\le t\le K\}$ is relatively compact in  $(C[0,K], \rho_K)$.
\end{proof}

\begin{lemma} \label{Lemmawk}
Let $M_k(t)=\max_{0\le s\le t} |w_k(s)|$.  Then   $w_k$ is absolutely continuous and,  for $t\notin \{\log (n/k): n\ge k\}$,
\be 
w_k'(t)= \frac{a}{c}(w_k(t)-w_k(t-L))+R_{k}(t),  \label{Equ1}
\ee 
where  the remainder term $R_k(t)$ satisfies
\be
|R_k(t)|\le C M_k(t). \label{boundRkt2}
\ee
Furthermore, there exists a sequence of intervals $\{E_k, k\ge 1\}$ in 
$[0, L+\log 5)$ with $|E_k|=O(1/k)$
such that,  for all $t\not \in E_k$,
\be
|R_k(t)|\le \frac{C}{ke^t} M_k(t). \label{boundRkt}
\ee
\end{lemma}

\begin{proof}
The absolute continuity of $w_k$ follows directly from its definition. For $t\in (\tau_{n-1,k}, \tau_{n,k})$,
\bestar
w_k'(t)=\frac{w_{n,k}-w_{n-1, k}}{\tau_{n,k}-\tau_{n-1,k}}&=&\frac{a}{m_n\log (n/(n-1))}(w_{k}(\tau_{n-1,k})-w_k(\tau_{n-m_n-1,k}))\\
&=& \frac{a}{c}(w_k(t)-w_k(t-L))+R_{k}(t),
\eestar
where 
\be 
R_k(t)=\frac{a}{m_n\log (n/(n-1))}(w_{k}(\tau_{n-1,k})-w_k(\tau_{n-m_n-1,k}))- \frac{a}{c}(w_k(t)-w_k(t-L)).  \label{Rkt}
\ee 
The bound  \eqref{boundRkt2}  follows immediately from  \eqref{Rkt}.

Finally, we prove \eqref{boundRkt}.
Define
\bestar
E_k=\left\{
\begin{array}{ll}
(L+\log (1-4/k), L+\log (1+4/k)), &  k>12/(1-c),\\\relax
[0, L+\log (1+4/k)), &  k\le 12/(1-c).
\end{array}
\right.
\eestar
Assume first that $k>12/(1-c)$ and
 $t\not \in E_k$.
If $t\in (\tau_{n-1,k}, \tau_{n,k})$,  then $n-1<ke^t<n$,  
$M_k(t)\ge \max_{n-1\ge \ell\ge k} |w_{\ell, k}|$,
and 
\be 
|\lfloor ke^{t-L}\rfloor-(n-m_n-1)|&=&|\lfloor ke^t(1-c)\rfloor-(n-m_n-1)|\nonumber\\
&\le& 2+|ke^t(1-c)-(n-cn-1)|\le 3+c<4. \label{diff}
\ee 
Since $t\not \in E_k=(L+\log (1-4/k), L+\log (1+4/k))$,  the integers $\lfloor ke^{t-L}\rfloor$ and $n-m_n-1$
 are either both greater than $k$ or both less than $k$.
Hence, by \eqref{windu},  \eqref{wkt} and \eqref{diff},
\bestar
&&|w_{k}(\tau_{n-1,k})-w_k(\tau_{n-m_n-1,k}) - (w_k(t)-w_k(t-L))|\\
&\le& |w_{k}(\tau_{n-1,k})- w_k(t)| +
|w_k(\tau_{n-m_n-1,k}) - w_k(t-L)|\\
&\le& 5\max_{n\ge \ell\ge n-m_n-5} |w_{\ell,k}-w_{\ell-1,k}|\\
&\le& \frac{C}{n-m_n-5} \max_{n-1\ge \ell\ge k} |w_{\ell, k}|
\le \frac{C}{n} M_k(t),
\eestar
where the last step uses the inequality $n-m_n-5>(1-c)n-6>(1-c)n/2$ 
since $n>k>12/(1-c)$.
Therefore,  for $t\in (\tau_{n-1,k}, \tau_{n,k})$,
\bestar
|R_k(t)|&\le& \Big|\frac{a}{m_n\log (n/(n-1))}-\frac{a}{c}\Big| |w_{k}(\tau_{n-1,k})-w_k(\tau_{n-m_n-1,k})|\\
&&+ \frac{|a|}{c}|w_{k}(\tau_{n-1,k})-w_k(\tau_{n-m_n-1,k}) - (w_k(t)-w_k(t-L))|\\
&\le& \frac{C}{n} M_k(t)\le  \frac{C}{ke^t} M_k(t),
\eestar
where we used $n-1<ke^t<n$, which holds since $t\in (\tau_{n-1,k}, \tau_{n,k})$,
and
\bestar
\frac{a}{m_n\log (n/(n-1))}-\frac{a}{c}=O(n^{-1}).
\eestar
Hence  \eqref{boundRkt} holds for $k>12/(1-c)$. The case $k\le 12/(1-c)$
is analogous. This completes the proof of Lemma \ref{Lemmawk}.
\end{proof}

Now we consider  the fundamental solution $w(t)$ of the equation \eqref{RDD0}, i.e.,  $w(t)$ satisfies
\be 
w'(t)=\frac{a}{c}(w(t)-w(t-L)),~~t\ge 0, \label{RDD1}
\ee 
together with the initial condition \eqref{initialcond}.
To study the properties of $w(t)$,   
we need to analyze its characteristic equation  
\be 
z=\frac{a}{c} (1-e^{-Lz}).  \label{chequa}
\ee

\begin{lemma} \label{lemma4}
The roots of \eqref{chequa} are
\be 
z_k=\frac{\lambda+W_k(-\lambda e^{-\lambda})}{L}, \qquad k\in \mathbb{Z},
\label{zk}
\ee 
where  $W_k$ denotes the $k$-th branch of the Lambert W function. 
Furthermore, the following hold:
\begin{itemize}
\item[$(a)$]  If $\lambda<0$, then  \eqref{chequa} has a unique real root $z_0=0$, and $\sup_{k\ne 0}\Re{z_k}<0$.
\item[$(b)$]  If $0<\lambda<1$, then \eqref{chequa} has  two real roots $z_0=0$ and $z_{-1}<0$, 
and $\sup_{k\ne 0}\Re{z_k}<0$.
\item[$(c)$]  If $\lambda=1$, then \eqref{chequa}  has a double real root at 
0, i.e., $z_0=z_{-1}=0$,  and $\sup_{k\not \in\{-1, 0\}}\Re{z_k}<0$.
\item[$(d)$]   If $\lambda>1$, then \eqref{chequa} has two  real roots $ z_{-1}=0$ and
\be
z_0=\frac{\lambda+W_0(-\lambda e^{-\lambda})}{L}:=\rho>0, \label{rho}
\ee  and $\sup_{k\not \in\{-1, 0\}}\Re{z_k}<0$.
\end{itemize}

\end{lemma}

\begin{proof} Let $w=L(z-a/c)$. It follows from \eqref{chequa} that $(z-a/c)e^{Lz}=-a/c$, and consequently,
\be 
we^w=L(z-a/c)e^{Lz}e^{-aL/c}=-(aL/c) e^{-aL/c}=-\lambda e^{-\lambda}. 
\label{wew}
\ee 
The roots of the equation are $w_k=W_k(-\lambda e^{-\lambda}), k\in\mathbb{Z}$.
Hence the roots of \eqref{chequa} are
\bestar
z_k=\frac{w_k}{L}+\frac{a}{c}=\frac{\lambda+W_k(-\lambda e^{-\lambda})}{L}, \qquad k\in \mathbb{Z}.   
\eestar

We shall use the following  basic properties of  the Lambert W function
(see \cite{Cermak2026, Corless1996}):
\begin{enumerate}
\item[(1)] For any $x\in\mathbb{R}$,  $\{W_k(x):k\not\in \{-1, 0\}\}$
are not real. 
\item[(2)] $W_0(x)\in \mathbb{R}$ if and only if $x\ge -e^{-1}$,
and $W_{-1}(x)\in \mathbb{R}$ if and only if $x\in [-e^{-1},0)$.
\item[(3)] If $x\in (-e^{-1},0)$,  then $W_{-1}(x)<-1<W_0(x)<0$. 
\item[(4)] $W_0(-\lambda e^{-\lambda})=-\lambda$ for $\lambda \le 1$,
and $W_{-1}(-\lambda e^{-\lambda})=-\lambda$ for $\lambda \ge 1$.
\item[(5)] For $x>0$,  $\sup_{k\ne 0}\Re(W_k(x))<\Re (W_0(x))$.
\item[(6)] For $x<0$,
$\sup_{k\not \in \{-1, 0\}}\Re(W_k(x))<\Re (W_{-1}(x))\le \Re (W_0(x))$.
\end{enumerate}

Since
\bestar
-\lambda e^{-\lambda}\in \left\{
\begin{array}{ll}
[0,\infty), & \lambda\le 0,\\
(-e^{-1}, 0), & 0<\lambda<1,\\
\{-e^{-1}\}, & \lambda=1,\\
(-e^{-1}, 0), & \lambda>1,
\end{array}
\right.
\eestar
the desired conclusions follow directly from properties
(1)--(6). 
\end{proof}

\begin{lemma} \label{lemma3.5}
There exists $\eta_0>0$ such that for any $\eta\in (0,\eta_0)$,
\bestar
w(t)=\left\{
\begin{array}{ll}
(1-\lambda)^{-1}+O(e^{-\eta t}), & \lambda<1,\\
2t/L+2/3+O(e^{-\eta t}),  & \lambda=1,\\
c_0e^{\rho t}+O(e^{(\rho-\eta) t}), & \lambda>1,
\end{array}
\right.
\eestar
where $c_0=1/(1-\lambda (1-c)^{\rho})$ and $\rho$ is defined in \eqref{rho}.
\end{lemma}

\begin{proof}
Let  $\hat{w}(s)=\int_0^{\infty} e^{-st}w(t)dt$ be the Laplace  transform of $w$. 
Then the Laplace transforms of $w'(t)$ and $w(t-L)$ are
\bestar
\int_0^{\infty} e^{-st}w'(t)dt=\int_0^{\infty}e^{-st}dw(t)
=-w(0)+s\int_0^{\infty} e^{-st}w(t)dt=s\hat{w}(s)-1,
\eestar
and
\bestar
\int_0^{\infty} e^{-st}w(t-L)dt=\int_{L}^{\infty} e^{-st}w(t-L)dt=e^{-sL}\int_{0}^{\infty} e^{-st}w(t)dt
=e^{-sL}\hat{w}(s),
\eestar
respectively.  Therefore, it follows from \eqref{RDD1} that
$s\hat{w}(s)-1=(a/c)(1-e^{-sL})\hat{w}(s)$.  Consequently,
\bestar
\hat{w}(s)=\frac{1}{D(s)},
\eestar 
where $D(s)=s-a(1-e^{-Ls})/c$.
By Lemma \ref{lemma4}, the roots of $D(z)=0$ are $z_k, k\in \mathbb{Z}$,
as defined in \eqref{zk}.
 Then, applying  the inverse Laplace transform, we obtain
\be 
w(t)=\frac{1}{2\pi \mathbf{i}} \int_{b-\mathbf{i}\infty}^{b+\mathbf{i}\infty} \frac{e^{st}}{D(s)}ds,
\label{wt}
\ee 
where $\mathbf{i}=\sqrt{-1}$ and $b$ is any real number satisfying $b>\sup_{k}\Re{z_k}$.

We only prove the case $\lambda=1$, the others being analogous and hence omitted.

Suppose $\lambda=1$, and let $b>0$ be a constant. By Lemma \ref{lemma4}, we can choose $\delta>0$
such that $\sup_{k\not \in \{-1,0\}} z_k<-\delta$ and $D(s)$ has only one double zero $s=0$ in the half-plane $\Re s>-\delta$. 
Let $\gamma_{T}$ denote the rectangular contour
with vertices $b-\mathbf{i}T, b+\mathbf{i}T, -\delta+\mathbf{i}T, -\delta-\mathbf{i}T$ and the positive counterclockwise orientation. Then we have
\bestar
 \frac{1}{2\pi \mathbf{i}}\int_{\gamma_T}\frac{e^{st}}{D(s)}\,ds=\frac{1}{2\pi \mathbf{i}}\int_{b-\mathbf{i}T}^{b+\mathbf{i}T}\frac{e^{st}}{D(s)}\,ds+R_1(T)+R_2(T),
\eestar
where
\bestar
 R_1(T):=-\frac{1}{2\pi }\int_{-T}^{T}  \frac{e^{(-\delta+\mathbf{i}y)t}}{D(-\delta+\mathbf{i}y)}\,dy,
\eestar
and
\bestar 
R_2(T):=\frac{1}{2\pi \mathbf{i}}\int_{-\delta}^b \Big(\frac{e^{(x-\mathbf{i}T)t}}{D(x-\mathbf{i}T)}
+\frac{e^{(x+\mathbf{i}T)t}}{D(x+\mathbf{i}T)}\Big)dx.
\eestar

For any $x, y\in \mathbb{R}$, we have
\be 
|D(x+\mathbf{i}y)|\ge |x+\mathbf{i}y-a/c|-|ae^{-L(x+\mathbf{i}y)}/c| \ge |y|-|ae^{-Lx}/c|.  \label{Dxy}
\ee 
This implies that  $|D(x+\mathbf{i}T)|, |D(x-\mathbf{i}T)|\to \infty$ uniformly for $x\in [-\delta, b]$,
and consequently,
\bestar
|R_2(T)|\le \frac{1}{2\pi}\int_{-\delta}^b \Big(\frac{e^{xt}}{|D(x-\mathbf{i}T)|}
+\frac{e^{xt}}{|D(x+\mathbf{i}T)|}\Big)dx \rightarrow 0,\qquad T\rightarrow \infty.
\eestar

By integration by parts, we obtain
\bestar 
R_{1}(\infty):=\lim_{T\rightarrow \infty}R_{1}(T)
= -\frac{e^{-\delta t}}{2\pi}\int_{-\infty}^{\infty}  \frac{e^{\mathbf{i}yt}}{D(-\delta+\mathbf{i}y)}\,dy
=\frac{e^{-\delta t}}{2\pi t}\int_{-\infty}^{\infty}  \frac{e^{\mathbf{i}yt}D'(-\delta+\mathbf{i}y)}{D^2(-\delta+\mathbf{i}y)}\,dy.
\eestar 
For $|y|>T_0:=|2ae^{L\delta}/c|$, it follows from \eqref{Dxy} that
$
|D(-\delta+\mathbf{i}y)|\ge |y|-|ae^{L\delta}/c|> |y|/2.
$
Note that $|D'(-\delta+\mathbf{i}y)|=|1-\lambda e^{L(\delta-\mathbf{i}y)}|\le 1+e^{L\delta}$
since $\lambda=1$.
Then  
\bestar
|R_{1}(\infty)|&\le& \frac{e^{-\delta t}}{2\pi t}\int_{-\infty}^{\infty}  \frac{|D'(-\delta+\mathbf{i}y)|}{|D^2(-\delta+\mathbf{i}y)|}\,dy \\
&\le& \frac{e^{-\delta t}( 1+e^{L\delta})}{2\pi t} \Big(\int_{|y|>T_0} \frac{4}{y^2}dy+\frac{2T_0}{M_0^2}\Big)=O(e^{-\delta t}/t), 
\eestar
where $M_0:=\inf_{|y|\le T_0}|D(-\delta+\mathbf{i}y)|>0 $.
Hence
\be
\lim_{T\to \infty}\frac{1}{2\pi i}\int_{\gamma_T}\frac{e^{st}}{D(s)}\,ds
=\int_{b-\mathbf{i}\infty}^{b+\mathbf{i}\infty}\frac{e^{st}}{D(s)}\,ds+O(e^{-\delta t}/t).  \label{Cont1}
\ee

By the residue formula,  
\bestar
\frac{1}{2\pi \mathbf{i}}\int_{\gamma_T}\frac{e^{st}}{D(s)}\,ds
=\operatorname{Res}_{s=0} \frac{e^{st}}{D(s)} 
=\lim_{s\to 0}\frac{d}{ds}\frac{s^2e^{st}}{D(s)}
=\frac{2t}{L}+\frac{2}{3}.
\eestar
Combining this with \eqref{wt} and \eqref{Cont1} yields that
for any $\eta\in (0,\delta)$
\bestar 
w(t)=\frac{2t}{L}+\frac{2}{3}+O(e^{-\eta t}).
\eestar 

The proof of Lemma \ref{lemma3.5} is complete.
\end{proof}

\begin{lemma} \label{lemma3.6}
For any $K\ge 1$,  
\be
\sup_{0\le t\le K} |w_{\lfloor ke^t\rfloor, k}-w(t)|\rightarrow 0,
\qquad k\rightarrow \infty.  \label{unifdist}
\ee
\end{lemma}

\begin{proof}
By the variation-of-constants formula (see \cite[Theorem 1.6.1]{Hale1993}),
\be 
w_k(t)=w(t)+\int_0^t w(t-s)R_k(s)ds.  \label{wktwt}
\ee 
Let  $M_k(t)$ be defined as in 
Lemma \ref{Lemmawk} and define $M(t):=\sup_{0\le s\le t}|w(s)|$.
It  follows from the proof of Lemma \ref{lemma3.2} that 
$\sup_k M_k(K)<\infty$.  Moreover,  $M(K)<\infty$ by the continuity of $w$.
Using \eqref{boundRkt2} and \eqref{boundRkt}, we obtain
\bestar
\sup_{0\le t\le K}|w_k(t)-w(t)|&\le& \sup_{0\le t\le K}\int_0^t |w(t-s)||R_k(s)|ds\\
&\le& CM(K) \Big( \int_0^K \frac{1}{ke^s} M_k(s)ds+\int_{[0,K]\cap E_k} M_k(s)ds\Big)\\
&\le&  CM(K) M_k(K)\Big( \int_0^{\infty} \frac{1}{ke^s} ds+|E_k|\Big)=O(1/k),
\eestar
where $E_k$ is the interval from  Lemma \ref{Lemmawk} and  $|E_k|=O(1/k)$. 
Hence, by the uniform continuity of $w$ on $[0,K]$, 
\bestar
\sup_{0\le t\le K} |w_{\lfloor ke^t\rfloor, k}-w(t)|&=&\sup_{0\le t\le K} |w_{k}(\log (\lfloor ke^t\rfloor/k))-w(t)|\\
&\le &\sup_{0\le t\le K}|w_k(t)-w(t)|+\sup_{0\le t\le K} |w(\log (\lfloor ke^t\rfloor/k))-w(t)|\rightarrow 0.
\eestar
This proves \eqref{unifdist}.
\end{proof}

\begin{lemma} \label{lemma3.7}
There exists $C>0$ such that
for $1\le k\le n$,
\be 
w_{n,k}\le \left\{
\begin{array}{ll}
C, & \lambda<1;\\
C(1+\log (n/k)), & \lambda=1;\\
C(n/k)^{\rho}, & \lambda>1,
\end{array}
\right.
\label{boundwnk}
\ee 
where $\rho$ is defined in \eqref{rho}.
Furthermore,  if $\lambda>1$, then there exist  $A_{k}\in \mathbb{R}$ and $\eta_0>0$ such that, for $n\ge k$ and $\eta\in (0,\eta_0)$,
\be 
w_{n,k}=A_{k}(n/k)^{\rho}+O((n/k)^{\rho-\eta}), \qquad A_{k}=\frac{1}{1-\lambda(1-c)^{\rho}}+O\Big(\frac{1}{k}\Big).  \label{boundwnk2}
\ee 

 \end{lemma}

\begin{proof}
We first prove \eqref{boundwnk}. For $\lambda>1$,
define $W_k(t):=\sup_{0\le s\le t} e^{-\rho s} |w_k(s)|$.
Then it follows from \eqref{wktwt} and Lemma \ref{lemma3.5} that, for any $s\in [0,t]$,
\be 
e^{-\rho s} |w_k(s)| &\le& e^{-\rho s}|w(s)|+\int_0^s |e^{- \rho (s-r)}w(s-r)|
 \cdot |e^{- \rho r}R_k(r)|dr\nonumber\\
&\le& C+C \int_0^s e^{- \rho r}|R_k(r)|dr. \label{Wkt}
\ee 
Let  $M_k(t)$ be defined as in 
Lemma \ref{Lemmawk}. Then $e^{- \rho r} M_k(r)\le \sup_{0\le \ell \le r} e^{- \rho \ell}|w_k(\ell)|
\le W_k(r)$. Hence, by Lemma \ref{Lemmawk},
\bestar
\int_0^s e^{- \rho r}|R_k(r)|dr &\le&  \int_0^s \frac{C}{ke^{r}} e^{-\rho r}M_k(r)dr
+\int_{r\in [0,s]\cap E_k} Ce^{- \rho r} M_k(r)dr\\
&\le & CW_k(s) \Big(\int_0^{\infty} \frac{1}{ke^{r}} dr
+ |E_k|\Big)\le \frac{C}{k} W_k(s).
\eestar
Substituting this into \eqref{Wkt} yields
\bestar
W_k(t)\le C\Big(1+\frac{1}{k}\Big) W_k(t).
\eestar
Hence,    $W_k(t) \le C$  holds for sufficiently large $k$.
This implies that $|w_k(s)|\le Ce^{\rho s}$, and hence
$w_{n,k}\le C(n/k)^{\rho}$   for all sufficiently large $k$.
The remaining finitely many values of $k$ can be handled using the inequality  $
 w_{n, k}\le  w_{n, j}w_{j,k},~k \le j \le n$, established in Lemma \ref{lemma3.1}. 
The cases  $\lambda<1$ and $\lambda=1$ can be treated similarly. More precisely, when $\lambda< 1$,  we set $W_k(t)=\sup_{0\le s\le t}  |w_k(s)|$  and use the boundedness of $w(t)$.  When $\lambda= 1$,  we instead set $W_k(t)=\sup_{0\le s\le t}  |w_k(s)|/(1+s)$  and use the linear growth of $w(t)$.  In both cases, the desired estimate \eqref{boundwnk} follows.

We next prove \eqref{boundwnk2}.  Set  $\tau_{n,k}=\log (n/k)$. Then, by the definition of 
$w_k$ and \eqref{wktwt}, we obtain
\be
w_{n,k}=w_k(\tau_{n,k})=w(\tau_{n,k})+ \int_0^{\tau_{n,k}} w(\tau_{n,k}-s)R_k(s)ds. \label{wnk}
\ee
We use the asymptotic expansion of $w(t)$ in Lemma \ref{lemma3.5}.
Assume that $\lambda>1$ and define  $R(t):=w(t)-c_0e^{\rho t}$, where
$c_0=1/(1-\lambda (1-c)^{\rho})$. 
Let $\eta_0$ be defined as in Lemma \ref{lemma3.5}, and let $\eta\in (0,\eta_0)$.
Then $R(t)=O(e^{(\rho-\eta) t})$.
It follows from \eqref{wnk} that
\bestar
w_{n,k}
&=& c_0\Big(\frac{n}{k}\Big)^{\rho}\Big(1+\int_0^{\tau_{n,k}} e^{-\rho s}R_k(s)ds\Big)+R(\tau_{n,k}) +\int_0^{\tau_{n,k}} R(\tau_{n,k}-s)R_k(s)ds\\
&=& A_{k}\Big(\frac{n}{k}\Big)^{\rho}+R_{n,k},
\eestar
where 
\bestar
&&A_{k}=c_0+c_0\int_0^{\infty} e^{-\rho s}R_k(s)ds,\\
&&R_{n,k}=R(\tau_{n,k}) +\int_0^{\tau_{n,k}} R(\tau_{n,k}-s)R_k(s)ds-c_0\Big(\frac{n}{k}\Big)^{\rho}\int_{\tau_{n,k}}^{\infty} e^{-\rho s}R_k(s)ds.
\eestar
By \eqref{boundwnk}, we have
$w_{n,k}\le C(n/k)^{\rho}$, and consequently
\bestar
|w_k(s)|\le \max\{w_{\lfloor ke^s\rfloor,k},  w_{\lfloor ke^s\rfloor+1,k})
\le C\Big(\frac{ke^s+1}{k}\Big)^{\rho}\le Ce^{\rho s}.
\eestar
This together with Lemma \ref{Lemmawk} implies that  $|R_k(s)|\le Ce^{\rho s}$ for $s\in E_k$, and
 $|R_k(s)|\le Ck^{-1}e^{\rho s-s}$  for $s\notin E_k$, where $E_k$ is
a interval in $[0, L+\log 5)$ satisfying $|E_k|=O(1/k)$. Hence
\bestar
\int_0^{\infty} e^{-\rho s}|R_k(s)|ds\le  C\Big(\int_0^{\infty} \frac{1}{ke^s}ds+\int_{E_k} ds\Big)=O\Big(\frac{1}{k}\Big).
\eestar
This yields $A_{k}=c_0+O(1/k)$. 

In what follows, we may assume without loss of generality that $0<\eta<1$. 
Then
\bestar
\int_0^{\tau_{n,k}} e^{(\rho-\eta) (\tau_{n,k}-s)}  |R_k(s)|ds
&\le& \Big(\frac{n}{k}\Big)^{\rho-\eta} \int_0^{\tau_{n,k}} e^{-(\rho-\eta) s}  |R_k(s)|ds \\
&\le&  C\Big(\frac{n}{k}\Big)^{\rho-\eta} \Big(\int_0^{\infty} \frac{1}{ke^{(1-\eta)s}}ds+\int_{0}^{L+\log 5} e^{\eta s}ds\Big)\le C\Big(\frac{n}{k}\Big)^{\rho-\eta},
\eestar
and
\bestar
 \int_{\tau_{n,k}}^{\infty} e^{-\rho s}|R_k(s)|ds
&\le& C\Big(\int_{\tau_{n,k}}^{\infty} \frac{1}{ke^s}ds+\int_{E_k \cap (\tau_{n,k},\infty)} ds\Big)\le   \frac{Ce^{-\tau_{n,k}}}{k}+CI(n\le  5e^Lk)\\
&\le &  C\Big(\frac{n}{k}\Big)^{-\eta}+C I(n\le  5e^Lk),
\eestar
where we used the fact that  $E_k \cap (\tau_{n,k},\infty)=\emptyset$
for $n> 5e^Lk$, and $E_k\subset [0, L+\log 5)$.
Therefore, by recalling that  $R(t)=O(e^{(\rho-\eta) t})$,
\bestar
|R_{n,k}|&\le&  Ce^{(\rho-\eta) \tau_{n,k}} 
 +C\int_0^{\tau_{n,k}} e^{(\rho-\eta) (\tau_{n,k}-s)}  |R_k(s)|ds+C\Big(\frac{n}{k}\Big)^{\rho}\int_{\tau_{n,k}}^{\infty} e^{-\rho s}|R_k(s)|ds\\
&\le& C\Big(\Big(\frac{n}{k}\Big)^{\rho-\eta}+1\Big)\le  C\Big(\frac{n}{k}\Big)^{\rho-\eta}.
\eestar
Hence
\bestar 
w_{n,k}=A_{k}(n/k)^{\rho}+O((n/k)^{\rho-\eta}).
\eestar
This proves  \eqref{boundwnk2}.
\end{proof}

\begin{lemma} \label{lemma3.7a}
Let   $\rho$ be defined as in \eqref{rho} when $\lambda>1$. 
\begin{enumerate}
\item[$(1)$] $p<p_c$  if and only if one of the following holds: (a) $\lambda\le 1$; (b) $\lambda>1$ and $\rho<1/2$. 
\item[$(2)$] $p=p_c$  if and only if $\lambda>1$ and $\rho=1/2$.
\item[$(3)$] $p>p_c$  if and only if $\lambda>1$ and $\rho\in (1/2,1)$.
\end{enumerate}
\end{lemma}

\begin{proof}
We first show that  $p\ge p_c$ implies  $\lambda>1$.

Indeed, if $p\ge p_c$, then
\bestar
\lambda=\frac{-(2p-1)\log (1-c)}{c}
\ge \frac{-(1+\sqrt{1-c})\log (1-c)}{2c}=\frac{-\log (1-c)}{2(1-\sqrt{1-c})}.
\eestar
Applying the inequality $-\log x>1-x$ for $x\in (0,1)$ with $x=\sqrt{1-c}$, 
we obtain 
$ -\log (1-c)>2(1-\sqrt{1-c})$, and consequently  $\lambda>1$.

We now consider the case $\lambda>1$. 
Define
\bestar
f(z)=\frac{a}{c} (1-e^{-Lz})-z,\qquad z\in [0,\infty).
\eestar
Then $f(z)$ is strictly concave on $(0,\infty)$ since $f''(z)=-\lambda Le^{-Lz}<0$.
By  Lemma \ref{lemma4},  $\rho$ is the unique positive solution of $f(z)=0$.
 Since $f(0)=f(\rho)=0$ and $f(z)$ is  strictly concave,
we have  $f(z)>0$ for  $z\in (0,\rho)$, and  $f(z)<0$ for  $z\in (\rho,\infty)$.
 In particular, $f(1/2)>0$ iff $\rho>1/2$.

A direct calculation gives
\bestar 
f(1/2)=\frac{a}{c} (1-e^{-L/2})-\frac{1}{2}=\frac{2p-1}{c}(1-\sqrt{1-c})-\frac{1}{2}
=\frac{2p-1}{1+\sqrt{1-c}}-\frac{1}{2}.
\eestar
Hence $f(1/2)=0$ is equivalent to 
\bestar 
\frac{2p-1}{1+\sqrt{1-c}}=\frac{1}{2},
\eestar
which is exactly $p=p_c$.  This proves (2).

Moreover,  $p>p_c$  if and only if $f(1/2)>0$, which is equivalent to  $\rho>1/2$. Since
\bestar
f(1)=\frac{a}{c} (1-e^{-L})-1=2p-2<0,
\eestar
we always  have $\rho<1$ by the concavity of $f$  and the fact that $f(0)=f(\rho)=0$.
Therefore, $p>p_c$  if and only if $\lambda>1$ and $\rho\in (1/2,1)$,
proving (3).

Finally,  (1) follows immediately from (2) and (3), and the fact that $\rho<1$
whenever $\lambda>1$.
\end{proof}

\begin{lemma}  \label{lemmasumw}
Let  $\sigma^2(c,p)$ and $\sigma_c^2$ be as defined in Theorem \ref{th1}.
Define
$B_n^2:=\sum_{k=1}^n w_{n,k}^2$.  

\begin{enumerate}
\item[$(1)$] If $p<p_c$, then  $B_n^2/n\rightarrow  \sigma^2(c,p)$.
 
\item[$(2)$] If $p=p_c$, then    $B_n^2/(n\log n)\rightarrow  \sigma^2_c$.

\item[$(3)$] If $p>p_c$, then
\be 
\frac{1}{n^{2\rho}}\sum_{k=1}^n w_{n,k}^2
\rightarrow \sum_{k=1}^{\infty} \frac{A_{k}^2}{k^{2\rho}},
\label{Bn03}
\ee 
where $\rho\in (1/2,1)$ and $A_{k}$ are defined in \eqref{rho} and Lemma \ref{lemma3.7}, respectively.
\end{enumerate}
\end{lemma}

\begin{proof}
We first prove (1). For any $\varepsilon\in (0,1)$, write
\be 
B_n^2=\sum_{k=\lfloor\varepsilon n\rfloor+1}^n w_{n,k}^2
+\sum_{k=1}^{\lfloor\varepsilon n\rfloor} w_{n,k}^2.
\label{Bn1}
\ee 

For fixed $\varepsilon\in (0,1)$, it follows from Lemma \ref{lemma3.6} that,
as $n\to \infty$,
\bestar
\max_{\varepsilon n<k\le n}|w_{n,k}-w(\log (n/k))|
\le \max_{k>\varepsilon n}\sup_{0\le t\le -\log \varepsilon} |w_{\lfloor ke^t\rfloor, k}-w(t)|\rightarrow 0.
\eestar
Hence
\be 
\frac{1}{n}\sum_{k=\lfloor\varepsilon n\rfloor+1}^n w_{n,k}^2
&=&\frac{1}{n} \sum_{k=\lfloor\varepsilon n\rfloor+1}^n w^2(\log (n/k)) +o(1)
=\int_{\varepsilon}^1 w^2(-\log t)dt+o(1)\nonumber\\
&=&\int_{0}^{-\log \varepsilon} \frac{w^2(u)}{e^u} du+o(1).  \label{Bn2}
\ee 
By Lemma \ref{lemma3.7} and  \ref{lemma3.7a},  if $p<p_c$, then
\be 
\frac{1}{n}\sum_{k=1}^{\lfloor\varepsilon n\rfloor} w_{n,k}^2\le g(\varepsilon),  \label{0toep}
\ee 
where $g(x)$ is   continuous  on $[0,1)$ and $g(0)=0.$
Indeed,  if $\lambda<1$,  then
\bestar
\frac{1}{n}\sum_{k=1}^{\lfloor\varepsilon n\rfloor} w_{n,k}^2\le C\varepsilon.
\eestar
If $\lambda=1$, then
\bestar
\frac{1}{n}\sum_{k=1}^{\lfloor\varepsilon n\rfloor} w_{n,k}^2
&\le& \frac{C}{n}\sum_{k=1}^{\lfloor\varepsilon n\rfloor} (1+\log (n/k))^2
\le \frac{C}{n}\sum_{k=1}^{\lfloor\varepsilon n\rfloor} \int_{k-1}^{k} (1+\log (n/x))^2 dx \\
&\le& \frac{C}{n}\int_{0}^{\varepsilon n} (1+\log (n/x))^2 dx 
= C\int_0^{\varepsilon} (1-\log u)^2du.
\eestar
If $\lambda>1$ and $\rho<1/2$, then
\bestar
\frac{1}{n}\sum_{k=1}^{\lfloor\varepsilon n\rfloor} w_{n,k}^2
&\le& \frac{C}{n}\sum_{k=1}^{\lfloor\varepsilon n\rfloor} (n/k)^{2\rho}
\le \frac{C}{n}\sum_{k=1}^{\lfloor\varepsilon n\rfloor} \int_{k-1}^{k} (n/x)^{2\rho} dx \\
&\le& \frac{C}{n}\int_{0}^{\varepsilon n} (n/x)^{2\rho} dx 
= C\int_0^{\varepsilon}u^{-2\rho} du.
\eestar
Hence \eqref{0toep} holds.  Letting $\varepsilon\to0$ in \eqref{Bn2} and \eqref{0toep} yields  (1). 

Next we prove (2) and (3). Assume $p\ge p_c$. Then  Lemma  \ref{lemma3.7a} gives $\lambda>1$ and $\rho\ge 1/2$. 
Moreover,  by Lemma \ref{lemma3.7},  there exist  $A_{k}$ and $\eta_0>0$ such that, for $n\ge k$ and $\eta\in (0,\eta_0)$,
\bestar
w_{n,k}^2=A_{k}^2(n/k)^{2\rho}+O((n/k)^{2\rho-\eta}), \qquad A_{k}=\frac{1}{1-\lambda(1-c)^{\rho}}+O\Big(\frac{1}{k}\Big).
\eestar
If $p=p_c$,  then $\rho=1/2$. In this case, we have
\bestar
\frac{1}{n\log n}\sum_{k=1}^n w_{n,k}^2&=&\frac{1}{\log n}\sum_{k=1}^n \frac{A_{k}^2}{k} + \frac{O(1)}{n^{\eta}\log n}\sum_{k=1}^n  k^{\eta-1}\\
&=& \frac{1}{(1-\lambda \sqrt{1-c})^2}\frac{1}{\log n}\sum_{k=1}^n \frac{1}{k}+\frac{O(1)}{\log n}\sum_{k=1}^n \frac{1}{k^2} +\frac{O(1)}{n^{\eta}\log n}\sum_{k=1}^n  k^{\eta-1}\\
&\rightarrow& \frac{1}{(1-\lambda \sqrt{1-c})^2}=\sigma_c^2.
\eestar
If $p>p_c$, then $\rho\in (1/2,1)$. Choosing $0<\eta<\min\{\eta_0, 2\rho-1\}$,  
we obtain
\bestar
\frac{1}{n^{2\rho}}\sum_{k=1}^n w_{n,k}^2
= \sum_{k=1}^n \frac{A_{k}^2}{k^{2\rho}}+\frac{O(1)}{n^{\eta}}\sum_{k=1}^n  k^{\eta-2\rho}
\longrightarrow \sum_{k=1}^{\infty} \frac{A_{k}^2}{k^{2\rho}}.
\eestar
The proof of Lemma \ref{lemmasumw} is complete.
\end{proof}

\section{Proof of main result}  \label{sect4}

Before proving the main result, we  give a result concerning the almost sure  and $L^2$ convergences of $S_n$.

\begin{lemma} \label{th0}
For $0<p<1$,  we have
\bestar
\frac{S_n}{n}\longrightarrow 0 ~~\text{a.s. and in}~ L^2.
\eestar
\end{lemma}

\begin{proof}
By Lemma  \ref{lemma1},   $S_n=\sum_{k=1}^n w_{n,k} \ep_k$,  where 
$\{\ep_k\}$ are defined in \eqref{innovation} and  $\{w_{n,k}\}$ are the weights  introduced in Lemma \ref{lemma1}.  For each $n\ge 1$,   $\{w_{n,k}\ep_k, 1\le k\le n\}$ is a martingale difference sequence. Consequently,  $\mbox{Var} (S_n)=\sum_{k=1}^n w_{n,k}^2 \mathbb{E}
(\ep_k^2)$.  Since
$\mathbb{E}(\ep_k^2)\le \mathbb {E}(X_k^2)=1$, 
it follows from Lemma \ref{lemmasumw} that
\bestar
\frac{\mbox{Var} (S_n)}{n^2} \le \frac{\sum_{k=1}^n w_{n,k}^2}{n^2}\rightarrow 0.
\eestar
Hence $S_n/n\xrightarrow{L^2} 0$.

Since $|\ep_k|\le 2$, applying the Azuma-Hoeffding inequality (see, for instance, Corollary 2.20 in \cite{Wainwright2019}) yields
\bestar
\mathbb{P}(|S_n|\ge n\varepsilon) \le 2\exp\Big\{-\frac{n^2\varepsilon^2}{8\sum_{k=1}^n w_{n,k}^2}\Big\}.
\eestar
By Lemma \ref{lemmasumw},  we have
\bestar
\sum_{n=1}^{\infty}\mathbb{P}(|S_n|\ge n\varepsilon) \le 2\sum_{n=1}^{\infty}\exp\Big\{-\frac{n^2\varepsilon^2}{8\sum_{k=1}^n w_{n,k}^2}\Big\}<\infty.
\eestar
Applying the Borel-Cantelli lemma then gives $S_n/n\to 0$~a.s.
\end{proof}

We now turn to the proof of the main result.

\begin{proof}[Proof of Theorem \ref{th1}]
Set $B_n^2=\sum_{k=1}^n w_{n,k}^2$ and write
\bestar
X_{n,k}=\frac{w_{n,k}\ep_k}{B_n},\qquad 1\le k\le n,
\eestar
 where $\{\ep_k\}$ are defined in \eqref{innovation} and  $\{w_{n,k}\}$ are the weights  introduced in Lemma \ref{lemma1}. Then, for any $n\ge 1$, 
$\{X_{n,k}, 1\le k\le n\}$ is a martingale difference sequence. Moreover, 
\bestar
\frac{S_n}{B_n}=\sum_{k=1}^n X_{n,k}.
\eestar  
To establish the asymptotic normality for $p<p_c$ and $p=p_c$, we apply the martingale central limit theorem (see, for instance,  Corollary 3.1 in 
\cite{Hall1980}) to the array $\{X_{n,k}\}$.

Assume that $p\le p_c$. By $|\ep_k|\le 2$ and  \eqref{boundwnk},  we have
\bestar
\max_{k\le n}|X_{n,k}|\le 2\max_{k\le n}|w_{n,k}|\le \frac{Ca_n}{B_n},
\eestar
where
\bestar
a_n=\left\{
\begin{array}{ll}
1, & \lambda<1,\\
1+\log n, & \lambda=1,\\
n^{\rho}, & \lambda>1 ~\mbox{and}~ \rho\le 1/2.
\end{array}
\right.
\eestar
By Lemmas \ref{lemma3.7a} and \ref{lemmasumw}, we have ${a_n}/{B_n}\to 0$. Then for any  $\varepsilon>0$, we have
$|X_{n,k}|\le \varepsilon$ for all $k\le n$ whenever $n$ is sufficiently large, and consequently,
\begin{equation}\label{cltLin}
\sum_{k=1}^n
\mathbb{E}\big(X_{n,k}^2I(|X_{n,k}|>\varepsilon) | \mathscr{F}_{k-1}\big)
=0.
\end{equation}
Therefore, the conditional Lindeberg condition is satisfied.

Next, we verify the convergence of the conditional variances. We have
\be 
\sum_{k=1}^n\mathbb{E}(X_{n,k}^2 | \mathscr F_{k-1})
=\frac1{B_n^2}\sum_{k=1}^nw_{n,k}^2(1-\mu_k^2) 
= 1-\frac1{B_n^2}\sum_{k=1}^nw_{n,k}^2\mu_k^2.  \label{condvar}
\ee
where 
\bestar 
\mu_k:=\mathbb{E}(X_{k}|\mathscr{F}_{k-1})=\frac{a}{m_k}(S_{k-1}-S_{k-m_k-1}).   
\eestar 
Since $m_k\sim ck$,    
applying Lemma \ref{th0} yields, as $k\rightarrow \infty$,
\be 
\mu_k=\frac{aS_{k-1}}{m_k}-\frac{aS_{k-m_k-1}}{m_k} \rightarrow 0~~\mbox{a.s.}  \label{mukconv}
\ee 
By Toeplitz's lemma, we get
\bestar
\frac1{B_n^2}\sum_{k=1}^nw_{n,k}^2\mu_k^2 \rightarrow 0~~\mbox{a.s.}
\eestar
This, together with \eqref{condvar}, implies 
 $\sum_{k=1}^n\mathbb{E}(X_{n,k}^2 | \mathscr F_{k-1})
\to 1$ a.s.
Applying the martingale central limit theorem now gives
$$
\frac{S_n}{B_n}\xrightarrow[\quad]{d}  N(0,1).
$$
 Combining this with Lemma \ref{lemmasumw} and Slutsky's theorem proves
\eqref{CLTsubcritical} and \eqref{CLTcritical}.

It remains to prove \eqref{superconv}. Assume that $\rho>1/2$ and
define
\bestar
\widetilde{R}_{n,k}=\frac{w_{n,k}}{n^{\rho}}-\frac{A_{k}}{k^{\rho}}.
\eestar
Then $\widetilde{R}_{n,k}=O(n^{-\eta}k^{-\rho+\eta}) $ by \eqref{boundwnk2}. Moreover, 
\bestar
\frac{S_n}{n^{\rho}}=\frac{1}{n^{\rho}}\sum_{k=1}^n w_{n,k}\ep_k
=\sum_{k=1}^n \frac{A_{k}\ep_k}{k^{\rho}}+\sum_{k=1}^n \widetilde{R}_{n,k}\ep_k.
\eestar
Choosing $0<\eta<\min\{\eta_0, \rho-1/2\}$ and using similar arguments to those in the proof of Lemma \ref{th0}, we have
\bestar 
\mbox{Var} \Big(\sum_{k=1}^n \widetilde{R}_{n,k}\ep_k\Big) = \sum_{k=1}^n \widetilde{R}_{n,k}^2 \mathbb{E}
(\ep_k^2) \le \sum_{k=1}^n \widetilde{R}_{n,k}^2 
\le C n^{-2\eta}\sum_{k=1}^{\infty}  k^{-2\rho+2\eta} \to 0.
\eestar 
This implies  $\sum_{k=1}^n \widetilde{R}_{n,k}\ep_k\xrightarrow{L^2} 0$.
Moreover, as in the proof of Lemma \ref{th0}, applying the Azuma–Hoeffding inequality and the Borel–Cantelli lemma yields
 $\sum_{k=1}^n \widetilde{R}_{n,k}\ep_k\to 0$ a.s.
Therefore,
\be 
\frac{S_n}{n^{\rho}}-\sum_{k=1}^n \frac{A_{k}\ep_k}{k^{\rho}}\longrightarrow 0 ~~\text{a.s. and in}~ L^2. 
\label{superconv1}
\ee 
On the other hand, we have
\bestar
\mathbb{E}\Big(\sum_{k=1}^n \frac{A_{k}\ep_k}{k^{\rho}}\Big)^2
=\sum_{k=1}^n \frac{A_{k}^2\mathbb{E}(\ep_k^2)}{k^{2\rho}}
\le \sum_{k=1}^{\infty} \frac{A_{k}^2}{k^{2\rho}}<\infty.
\eestar
Hence, by the $L^2$ convergence theorem for martingales (see \cite[Theorem 4.4.6]{Durrett2019}),
\be 
\sum_{k=1}^n \frac{A_{k}\ep_k}{k^{\rho}} \longrightarrow W ~~\text{a.s. and in}~ L^2,   \label{superconv2}
\ee 
where $W=\sum_{k=1}^{\infty}  A_{k}\ep_k/k^{\rho}$ is a real-valued random variable.  Now \eqref{superconv} follows immediately from \eqref{superconv1}
and \eqref{superconv2}.
\end{proof}

\end{document}